\documentclass[final]{siamltex}

\usepackage {graphicx}
\title{On the problem of global optimisation of a multivariable function.}
\author{Michael M. Medynski \\
\it \scriptsize Moscow Aviation Institute (State Technical University), 
Russia 
\\medmm@hotmail.com}
\begin{document}

\maketitle

\begin{abstract}
One of the actual problems in the field of numerical optimisation, as is 
well known, is the problem of the search for the global extremum of a 
multivariate function [1-9,13,14,17-21]. 
Various versions of the random 
search methods [6,8,9] are considered to be the most reliable to solve the 
problem of global optimisation. In this work we present the little-known 
methods of Halton and LP-search, which has been proved as 
one of the best practical 
solutions of the global optimisation problem.
\end{abstract}

\begin{keywords}\end{keywords}
\begin{AMS}\end{AMS}
\pagestyle{myheadings}
\thispagestyle{plain}
\markboth{M. Medynski}{On the problem of global optimisation...}

\section{Introduction}

It is known [1 - 6, 13, 14, 21], that many problems in the engineering 
practice are related to the search for the best solution in the set of all 
allowable, and can be reduced to a problem of nonlinear programming:
$$
{\rm to\;\;\;find} \;\;\; {\rm \bf min\;\;  f} (\bar {x}),
\qquad \bar {x} \in R^{n }$$
given the restrictions in form of equalities and inequalities defining the 
so-called area of allowable solutions \textbf{G}:
$$
\left\{
\begin{array}{l}
{\rm \bf h}_{k}(\bar {x}) = 0,\quad k =1,2,\ldots ,m, \\
{\rm \bf g}_{k}(\bar {x})\le  0,\quad  k = m+1,m+2,\ldots ,p.
\end{array} 
\right.
$$

Objective function f($\bar {x})$= f($x_{1}$,$x_{2}$,\ldots ,$x_{n})$ 
defines the objective for choosing the best solution from a set of 
alternatives. (One can speak of minimization, since maximization of f($\bar 
{x})$ is obviously equivalent to minimization of function (- f($\bar 
{x}))$). Any vector $\bar {x}$, satisfying the restrictions, is called an 
allowable vector or an allowable point.

The allowable vector $\bar {x}^{\ast }$ = ($x_{1}^{\ast }$, 
$x_{2}^{\ast }$, \ldots , $x_{n}^{\ast })$, which imposes the 
minimum of the objective function f($\bar {x})$, is called the optimal 
point, and the appropriate value of function f($\bar {x}^{\ast })$ -- the 
optimal value of the objective function. The pair $\bar {x}$* and f($\bar 
{x}$*) is the optimal solution. A distinction is made between local and 
global optimum solutions. In both cases f($\bar {x}^{\ast }) 
\le $ f($\bar 
{x})$, but for the global optimal solution this condition is 
satisfied for all $\bar {x} \in $G, while for the local optimal solution the 
condition is satisfied only for a small vicinity of a point $\bar 
{x}^{\ast }$.

The global minimum is the optimal solution for the whole set of allowable 
solutions. It is better than the other solutions corresponding to local 
minima, and, as a rule, it is the one that is required to be found. However, 
all known and effectively working numerical methods of optimisation solving 
the above-mentioned problem of nonlinear programming are local, i.e. 
determining only locally optimal solutions [1 -- 9]. This fact proves to be 
true by the numerical experiments that have been carried out (including the 
algorithms used in the mathematical packages MATLAB 6 and MAPLE 6).

\section{Halton {\&} LP - search methods}

The methods of Halton and LP-search are the deterministic analogues of the 
global random search. It is known, that the elementary global search of a 
point of the minimum $\bar {x}$* of the objective function f ($\bar {x})$ is 
performed as follows: in the allowable area G some trial points $\bar 
{x}^{1}$,$\bar {x}^{2}$, \ldots , $\bar {x}^{N}$ are chosen, at each 
one the value of the objective function is calculated, and the point is 
selected at which the objective function attains the smallest value. This 
point serves as the first approximation of the required point $\bar {x}^*$ of 
the global minimum of the objective function f($\bar {x})$. It is 
considered, that the search method converges, if at least one of the points 
$\bar {x}^{1}$, $\bar {x}^{2}$, \ldots , $\bar {x}^{N}_{, }$given a 
big enough N, will fall into some predefined small vicinity U of the point 
of the global minimum of the objective function f($\bar {x})$.

Let the allowable area G to have a form of a singular n-dimensional cube 
K$^{n}$ ={\{}0$ \le x_{j} \le $1, j = 1,2,3,\ldots ,n{\}}. As the 
trial points $\bar {x}^{1}$, $\bar {x}^{2}$, \ldots , $\bar 
{x}^{N}_{ }$ for a random search we shall choose the independent random 
points evenly distributed in K$^{n}$. The probability that at least one 
point will fall in a small vicinity of a point of minimum U, equals 
P=1--(1--U)$^{N}$ and approaches 1 when N $ \to  \quad \infty $, i.e. the method 
converges.

For a non-random search non-random points $\bar {x}^{1}$, $\bar 
{x}^{2}$, \ldots ,$\bar {x}^{N}$ are chosen as the trial points. At the 
first glance it seems, that points should be chosen in regular intervals on 
area K$^{n}$ (Fig. 1). However it is not so [13, 14]. Suppose, that function 
f($\bar {x})$ depends only on one argument f($\bar {x})$ = f($x_1 )$.

\begin{figure}[t!]
\includegraphics[width=0.5\textwidth]{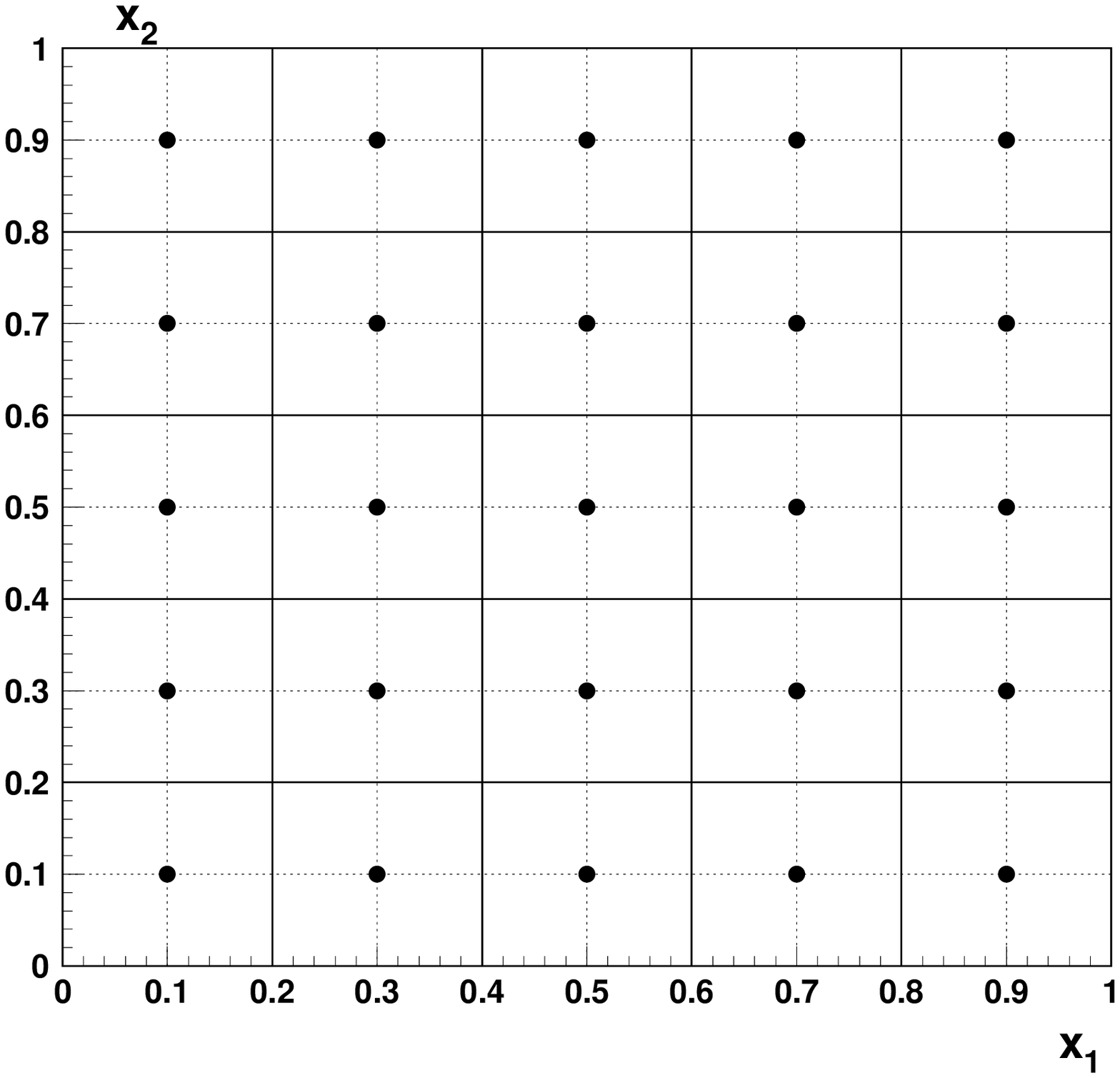}
\hfill
\includegraphics[width=0.5\textwidth]{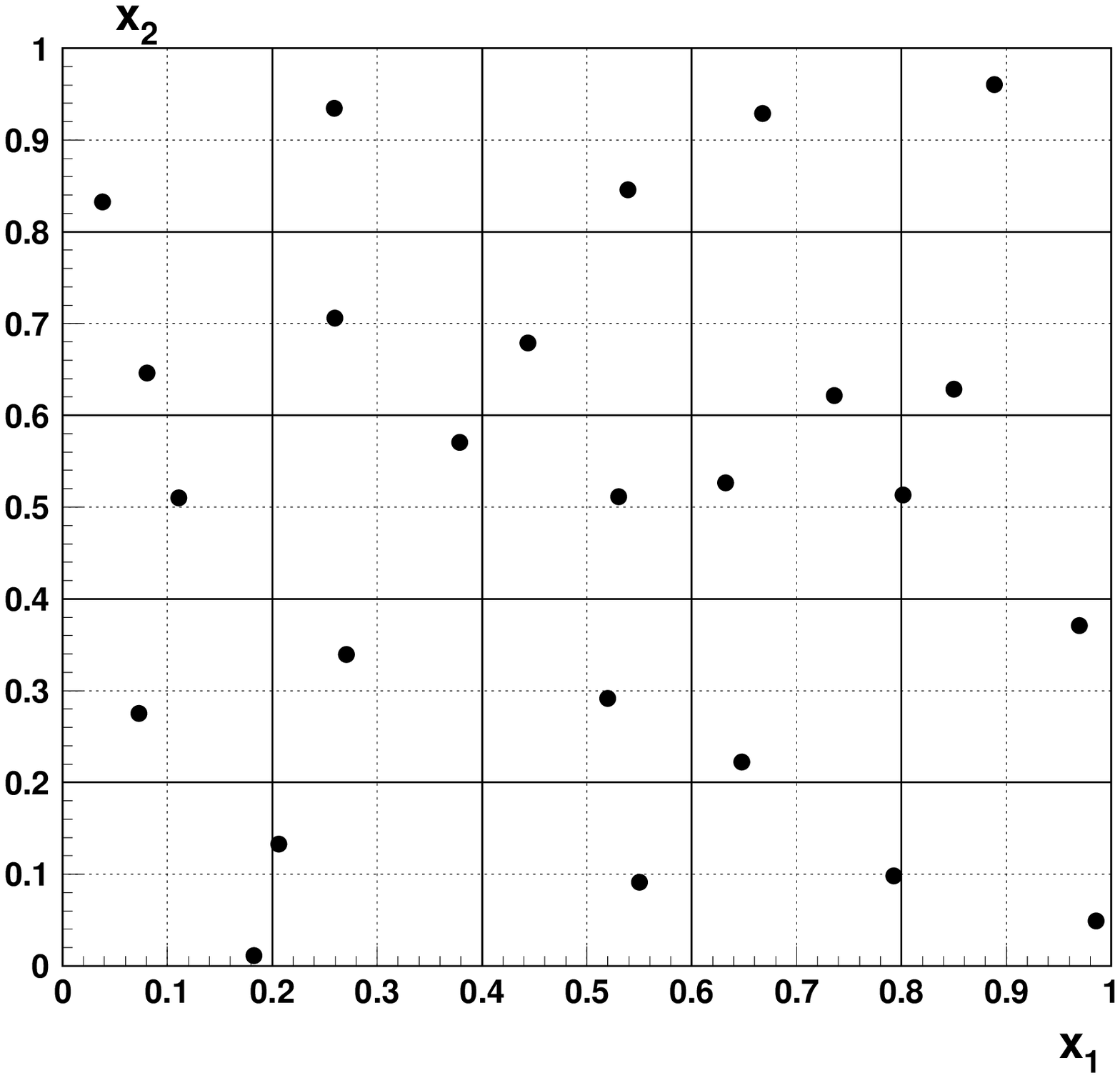}
\\
\parbox[t]{0.47\textwidth}{\caption{
"Evenly grid" for n=2 (N=25)
}\label{fig6}}
\hfill
\parbox[t]{0.47\textwidth}{\caption{Improved grid for 
n=2 (N=25)
}\label{fig7}}
\end{figure}

Calculating its values at points of such a "evently 
 grid" consisting of N = M$^{n}$ 
points (Fig. 1), we obtain M various values among which we shall choose the 
smallest, i.e. N - M points were calculated to no purpose. So, for 
example, 
for the case of a function of two variables f($x_1 ,x_2 )$, defined in 
K$^{2}$, which strongly depends only on $x_1 $, we shall obtain only five 
various values, each of which is repeated five times (Fig. 1). If f($\bar 
{x})$ "essentially" depends on $m < n$ arguments, N - M$^{m}$ of m points 
will 
be 
considered vainly. So, for example, if n = 4, m = 2 and M = 10 the number of 
"superfluous" values of function will be equal to N - M$^{m}$ = M$^{n}$ - 
M$^{m}$ = 10$^{4}$ - 10$^{2}$ = 9900 while the number of calculated values N 
= 10000! For the first time such feature was mentioned in 1957 in [10].

Random grids are devoid of the indicated flaw, for the probability of 
occurrence of identical abscissa at several random points is equal to zero. 
But can it be "good" to pick and choose the trial points, instead of relying 
on the will of a chance in choosing them? For example, if n = 2, points, 
represented in Fig. 2, are projected onto each of the coordinate axes to N 
various points, and, calculating the values of the function at the points of 
such a grid, we shall acquire N various values giving a much better 
representation of the function and of the range of its variation.

The idea of Halton's and LP-search methods is to use as trial points the 
points from the evenly distributed sequence of points in $K^{n}$ 
[11,14-16], 
namely: Halton's sequence $P_{1}$, $P_{2}$,\ldots , $P_{N}$,\ldots , 
constructed in [15, 16] and LP$_{\tau }$ -sequence $Q_{1,}Q_{2, }$\ldots 
,$Q_{N}$, \ldots , constructed in works of I.M.Sobol [11]. The Halton and 
LP$_{\tau }$-sequences are the deterministic analogues of a sequence of 
independent random points evenly distributed in K$^{n}$. Global search over 
the points of these sequences converges: it has been proved, that if S$_{N 
}$(U) -- is the number of the points which have fallen into a small vicinity 
U of the point of the minimum of the objective function f ($\bar {x})$, then 
S$_{N}$ (U) $ \ge $ 1 is satisfied, given a big enough N. The more evenly 
the points are distributed, the faster the search converges. It is 
considered, that the points of the LP$_{\tau }$-sequence are the most evenly 
distributed compared to all other known sequences [11]. (At the same 
time, test calculations, the results of which are presented below, have not 
shown any substantial advantage of the points of the LP$_{\tau }$-sequence 
over the Halton's sequence points).

\subsection{Halton's Sequence}

If r$_{1}$, r$_{2}$, \ldots , r$_{n }$- are pairs of coprimes, then Halton's 
sequence is defined as the sequence of points {\{}P$_{i}${\}} in K$^{n}$ 
with Cartesian coordinates P$_{i}$ = (p$_{r_1} $( i ), p$_{r_2} $( i ), 
\ldots , p$_{r_n} $( i )), i = 1, 2, \ldots , where p$_{r }$( i ) is the 
numerical sequence determined as follows.

If a$_{m}$a$_{m - 1}$\ldots a$_{2}$a$_{1 }$ is an integer i written in base 
r, then p$_{r }( i ) = 0.a_{1}a_{2}\ldots a_{m - 1}a_{m}$ (here all 
the $a_{s}$ are whole ``r'' numbers, i.e. are equal to one of the values 0, 
1, \ldots , r - 1). Accordingly, in the decimal system we have
$$
i = \sum\limits_{s = 1}^m {a_s } r^{s - 1}; \quad p_{r} ( i ) = 
\sum\limits_{s = 1}^m {a_s } r^{ - s}$$
For example, the first 10 values $p_{2} ( i )$, $p_{3} ( i )$, $p_{5} ( i 
)$, $p_{7} ( i )$, $p_{11} ( i )$ are equal:

\begin{table}[h]
\begin{tabular}
{|c|c|c|c|c|c|c|c|c|c|c|}
\hline
i& 
1& 
2& 
3& 
4& 
5& 
6& 
7& 
8& 
9& 
10 \\
\hline
P$_{2}$ ( i )& 
1/2& 
1/4& 
3/4& 
1/8& 
5/8& 
3/8& 
7/8& 
1/16& 
9/16& 
5/16 \\
\hline
P$_{3}$ ( i )& 
1/3& 
2/3& 
1/9& 
4/9& 
7/9& 
2/9& 
5/9& 
8/9& 
1/27& 
10/27 \\
\hline
P$_{5}$ ( i )& 
1/5& 
2/5& 
3/5& 
4/5& 
1/25& 
6/25& 
11/25& 
16/25& 
21/25& 
2/25 \\
\hline
P$_{7}$ ( i )& 
1/7& 
2/7& 
3/7& 
4/7& 
5/7& 
6/7& 
1/49& 
8/49& 
15/49& 
22/49 \\
\hline
P$_{11}$ ( i )& 
1/11& 
2/11& 
3/11& 
4/11& 
5/11& 
6/11& 
7/11& 
8/11& 
9/11& 
10/11 \\
\hline
\end{tabular}
\label{tab1}
\end{table}

\noindent
All such sequences are evenly distributed in K$^{n}$. 
In practice, as a rule, the first \textbf{n} prime numbers are chosen as 
r$_{1}$, r$_{2}$, \ldots , r$_{ n}$: r$_{1 }$= 2, r$_{2}$ = 3, r$_{3}$ = 5, 
\ldots also used are n-dimensional points P$_{i}$ = ( p$_{2}$( i ), p$_{3}$( 
i ), \ldots ,p$_{rn}$( i )), i = 1,2,3,\ldots

For example:

\begin{enumerate}
\item
at n = 2, we obtain Halton's sequence as following:

P$_{1}$ = (1/2, 1/3), P$_{2}$ = (1/4, 2/3), P$_{3}$ = (3/4, 1/9), \ldots , 
P$_{10}$ = (5/16, 10/27), \ldots ;

\item
at n = 3, we obtain:

P$_{1}$ = (1/2, 1/3, 1/5), P$_{2}$ = (1/4, 2/3, 2/5), \ldots , P$_{10}$ = 
(5/16, 10/27, 2/25), \ldots ;

\item
at n = 4, we obtain:

P$_{1}$ = (1/2, 1/3, 1/5, 1/7), P$_{2}$ = (1/4, 2/3, 2/5, 2/7), \ldots , 
P$_{10}$ = (5/16,10/27,2/25,22/49), \ldots 

\end{enumerate}

\subsection{LP$\tau $ -sequence }

If in the binary system the index of a point is entered as: i = e$_{m}$e$_{m 
- 1}$\ldots e$_{2}$e$_{1}$, where e$_{1, }$e$_{2}$, \ldots , e$_{m}$ -- are 
binary numbers, each equal to 0 or 1, then the corresponding decomposition 
of the index of a point in the decimal numeration system has the following 
form: i = 2$^{m - 1 }$e$_{m }$+ 2$^{m - 2 }$e$_{m - 1 }$+ \ldots + 2e$_{2}$ 
$_{ }$+ e$_{1})$, for all 1 $ \le $ j $ \le $ n

$$
q_{j}^{i }= e_{1 }V_{j}^{1 } *  e_{2 }V_{j}^{2 } *  
e_{3 }V_{j}^{3 } * \ldots  *  e_{m }V_{j}^{m},
$$

\noindent
where q$_{j}^{i}$, j = 1,2,\ldots ,n -- are coordinates of points Q$_{i}$ 
of LP$_{\tau }$ - sequence.

Vectors $\bar {V}^{S}$ = ($V_{1}^{S}$ , $V_{2}^{S}$ , \ldots , 
$V_{n}^{S}$ ) are easily calculated.

In work [11] a table of numerators of the point coordinates of $\bar 
{V}^s$ for 1$ \le $ n $ \le $ 13, 1$ \le $ s $ \le $ 20 is presented, 
referred to as directing points. Denominators of all coordinates of a point 
$\bar {V}^s$ equal 2$^{s}$. For example, at n = 7, s = 3 we have $\bar 
{V}^3$= (1/8, 5/8, 7/8,1/8, 5/8,7/8,3/8). Thus, the number of points 
$N < 2^{21}$, 
dimension n $ \le $ 13. In work [14] a table of directing points 
for 1$ \le $ n $ \le $ 51, 1 $ \le $ s $ \le $ 20 is presented. 
Unfortunately, there are no simple-enough formulas for calculation of 
directing points for any \textbf{n} and \textbf{s}.

Operation * means digit-by-digit addition by the module 2 in the binary 
system:
$$
{\rm if}\;\; a = \sum\limits_{k = 1}^l a _{k }2 ^{ - k };\;\; b = \sum\limits_{k = 
1}^l b _{k }2 ^{ - k };\;\; {\rm than} \;\; a  *  b = \sum\limits_{k = 1}^l c _{k 
}2 ^{ - k}, \;\; {\rm where} \;\; c _{k } = (a _{k }+ b _{k }) {\rm mod}_2.
$$
(In more detail: to calculate the "sum" a * b, it is necessary to write down 
both these numbers in binary notation and to combine digits in the 
appropriate places by rules:
0 + 0 = 0, 1 + 0 = 0 + 1 = 1, 1 + 1 = 0, i.e. without carrying units to the 
higher digit places.
For example: $\frac{7}{8} \quad  *  \quad \frac{11}{16}$ = 0.111 $ * $ 0.1011 = 
0.0101 = $\frac{5}{16})$.

In all computers there is a special command for carrying out the operation 
*. The number of operations performed by a computer for the calculation of 
Q$_{i}$ increases with the increasing i, but slowly, as log$_{2}i$.

For example, we have to calculate points Q$_{13}$ and Q$_{22}$ in a 
4-dimensional cube.

Since in binary notation decimal number 13 is represented as 1101, then 
e$_{1}$ = 1, e$_{2}$ = 0, e$_{3}$ = 1, e$_{4}$ = 1 and coordinates of 
Q$_{13}$ equal: q$_{j}^{13}=V_{j}^{1} \quad  *  \quad V_{j}^{3} \quad  * $ 
$V_{j}^{4}$ , 1 $ \le $ j $ \le $ 4

Vectors $\bar {V}^{1}$ , $\bar {V}^{3}$ , $\bar {V}^{4}$ could be 
obtained with the help of the table: 

$\bar {V}^{1 \quad } = (\frac{1}{2}$, $\frac{1}{2}$, $\frac{1}{2}$, 
$\frac{1}{2})$; $\bar {V}^{3}$ = ($\frac{1}{8}$, $\frac{5}{8}$, 
$\frac{7}{8}$, $\frac{1}{8})$; $\bar {V}^{4}$ = ($\frac{1}{16}$, 
$\frac{15}{16}$, $\frac{11}{16}$, $\frac{5}{16})$.

\noindent For the coordinates in the binary system, we have

q$_{1}^{13}$ = 0.1 $ * $ 0.001 $ * $ 0.0001 = 0.1011 = $\frac{11}{16}$; 
q$_{2}^{13}$ = 0.1 $ * $ 0.101 $ * $ 0.1111 = 0.1101 = $\frac{13}{16}$ ;

q$_{3}^{13}$ = 0.1 $ * $ 0.111 $ * $ 0.1011 = 0.1101 = $\frac{13}{16}$; 
q$_{4}^{13}$ = 0.1 $ * $ 0.001 $ * $ 0.0101 = 0.1111 = $\frac{15}{16}$.

\noindent So, Q$_{13}$ = (11/16, 13/16, 13/16, 15/16).

Similarly, since in the binary system the number 22 can be written down as 
10110, then Q$_{22}=\bar {V}^2 * \bar {V}^3 * \bar {V}^5$, and the 
coordinates of Q$_{22}$ equal:

q$_{1}^{22}$ = 1/4$ * $1/8$ * $1/32 =0.01$ * $0.001$ * $0.00001 = 0.01101 
= 13/32,

q$_{2}^{22}$ = 3/4$ * $5/8$ * $17/32 = 0.11$ * $0.101$ * $0.10001 = 
0.11101 = 29/32,

q$_{3}^{22}$ = 1/4$ * $7/8$ * $13/32 = 0.01$ * $0.111$ * $0.01101 = 
0.11001 = 25/32,

q$_{4}^{22}$ = 3/4$ * $1/8$ * $31/32 = 0.11$ * $0.001$ * $0.11111 = 
0.00011 = 3/32.

\noindent
Thus, Q$_{22}$ = (13/32, 29/32, 25/32, 3/32).

Using points Q$_{i}$, it is possible to obtain the appropriate trial points 
$\bar {x}^{i}$ = ($x_{1}^{i}$ , $x_{2}^{i}$ , \ldots 
,$x_{n}^{i}$ ) in any parallelepiped $\Pi$: 
a $_{j}  \le  x_{j}  \le 
b_{j}$ , $1  \le  j  \le  n$. Their coordinates can be calculated by 
the formula

$$x_{j}^{i} = a_{j } + q_{j}^{i} (b_{j} - a_{j}) , \quad 1  \le  j 
 \le  n.$$

It has been proved in [14], that if points Q$_{i }$= ( q$_{1}^{i}$, 
q$_{2}^{i}$, \ldots , q$_{n}^{i}$ ) form an evenly distributed sequence 
in K$^{n}$, then points $\bar {x}^{i}$ = ($x_1 ^{i}$, $x_2 
^{i}$,\ldots , $x_n ^{i}$) also form an evenly distributed sequence in 
\textbf{n}-dimensional parallelepiped $\Pi$. Also it has been proved, that if 
$\bar {x}^1$, $\bar {x}^2$, \ldots , $\bar {x}^i$, \ldots is a sequence of 
points evenly distributed in a given parallelepiped $\Pi$ with the sides 
parallel to the coordinate axes, and G $ \subset \Pi$ -- is some area with 
the positive volume $V_{G }> 0$, then all selected points $\bar {x}^i \in 
 G$ , also form a sequence of points evenly distributed in $G$. The ratio of 
volumes $\gamma  = V_{G}/V_{\Pi}$ is called the efficiency of selection, 
as, on the average, to obtain one point in G it is necessary to consider 
1/$\gamma $ points in $\Pi$.

Naturally there is a question: why are these "quasi-random numbers" better 
than the random numbers in the Monte Carlo method [10,12]?

\begin{enumerate}
\item
On all the classes of test functions f($\bar {x})$ on which numerical 
experiment was carried out, Halton and LP-search methods allowed to attain 
the same accuracy as the random search, but with the number of trial points 
being 2-4 times less, and with faster convergence.
\item
Calculation of points of Halton and LP$_{\tau }$ - sequences is rather 
simple and can be easily programmed.
\end{enumerate}

\section{Test problems and estimation of efficiency of Halton and 
LP-search methods}

In practice the estimation of efficiency of numerical methods of 
optimisation is usually performed through computing experiments aimed at 
solving the so-called special test problems [2, 3]. One of the major 
questions here is the question on how widely can the method be used, i.e. 
whether it is possible to solve the majority of problems with the help of a 
given method?

Methods were checked using the following test functions:

\begin{enumerate}

\item
"Rozenbrock's Function":
$$
f (x_{1},x_{ 2}) = 100(x_{2} - x_{1}^{2})^{2} + (1 - 
x_{1} )^{2}$$ 
--- unimodal function of a "gully" type with nonlinear 
"bottom" of a parabolic kind (has an abrupt rounded hollow along a curve 
$x_{2 = }x_{1}^{2})$. The minimum of the function is reached at the point 
$\bar {x}^{\ast }$ = (1,1) and $f(\bar {x}^{\ast })$ = 0.

\item
"Fletcher-Powell's Function":
$$
f(x_{1} ,x_{2} ,x_{3}) = 100{\{}[x_{3} - 5\theta  
(x_{1},x_{2})]^{2} + [(x_{1}^{2}+x_{2}^{2})^{1 / 
2} - 1]^{2}{\}} + x_{3}^{2},$$
where 
$$
\theta (x_1 ,x_2 )=\left\{ {\begin{array}{l}
\displaystyle \frac{\mbox{1}}{\pi }\mbox{arctg}(\frac{x_2 }{x_1 }), 
\quad \mbox{ if }x_1 > 0, \\ 
\displaystyle  1 + \frac{1}{\pi }\mbox{arctg}(\frac{x_2 }{x_1 }), 
\quad \mbox{if  }x_1 < 0, \\ 
 \end{array}} \right.$$
- function with a "coiling ravine". Its partial derivatives of the 1st 
order are 
piecewise continuous. The global minimum is reached at the point $\bar 
{x}^{\ast }$ = (1,0,0) è f ($\bar {x}^{\ast })$ = 0.

\item
"Powell's Function":

$$
f (x_{1} ,x_{2} ,x_{3} ,x_{4}) = (x_{1} + 
10x_{2})^{2} + 5(x_{3} - x_{4})^{2} + (x_{2} - 
2x_{3})^{4} + 10(x_{1} - x_{4})^{4}
$$
- function has a "flat 
gully" bottom (weak singularity). The minimum is reached at the point $\bar 
{x}^{\ast }$ = (0,0,0,0) and f ($\bar {x}^{\ast })$ = 0. The Hessian of 
the function is singular at the point $\bar {x}^{\ast }$.

\item
"Wood's Function":

$$
f(x_{1},x_{2},x_{3},x_{4}) = 100(x_{2 }- 
x_{1}^{2})^{2} + (1 - x_{1})^{2 }+ 90(x_{4} - 
x_{3}^{2})^{2} + (1 - x_{3})^{2 }
+ 10.1((x_{2 }- 1)^{2} + (x_{4} - 1)^{2}) +
$$
$$
\qquad\qquad\qquad
+ 19.8(x_{2 }- 
1)(x_{4} - 1)
$$
- function has some local minima, distinct from the global one. The global 
minimum is achieved at the point $\overline x ^{\ast }$ = (1,1,1,1) and f 
($\overline x ^{\ast })$ = 0.

\item 
$$
f (x_{1} ,x_{2} ) = x_{1}^{2}+x_{2}^{2} - \cos 
(18x_{1}) - \cos (18x_{2} )$$ - the function has 25 local minima on $[- 
\pi , \pi ]$. The global minimum is achieved at a point $\bar {x}^{\ast 
}$ = (0,0) and its value is equal to $f(\bar {x}^{\ast } )$ = -2.

\item

"Himmelblau's Function": (10 variables)

$$
f(\bar {x})=\sum\limits_{i = 1}^{10} {((\ln (x_i } - 2))^2 + (\ln (10 - 
x_i ))^2) - (\prod\limits_{i = 1}^{10} {x_i } )^{0.2}$$

\noindent
with boundary conditions 2.001$ < x_i <$ 9.999, i = 1, 2, \ldots , 10 
(outside of allowable area objective function is undefined). The minimum is 
achieved at the point $\bar {x}^{\ast }$ = (9.351; 9.351; 9.351; 9.351; 
9.351; 9.351; 9.351; 9.351; 9.351; 9.351) and $f ( \bar {x}^{\ast })$= - 
45.778.

\item
$$
 f(x_1 ,x_2 ) =  \frac{1 + x_1 }{x_1 x_2 ^2} \biggl(25(1+x_1 )+0.5  \cdot 
$$ $$\sqrt {0.133 \cdot 10^7 + 
40931.68x_1 ^2 + 999.44x_2 ^4 - 32613.30x_2 ^2 + 12543.58x_1 x_2 ^2 - 
122795.04x_1 } \biggr)^{2},
$$
where  0.1 $< 
x_1 < $ 5; 0.1 $< 
x_2 <$ 10-- objective function from the practice of design of automated 
drives for flying devices (V.A. Polkovnikov Electric, Hydraulic and 
Pneumatic drives for flying devices and the range of their dynamic 
capabilities. M.: Publ. MAI, 2002).

\end{enumerate}

Any optimisation procedure should solve the specified problems effectively.

The choice of both multimodal, and unimodal test functions is not 
accidental. It is necessary to be convinced of efficiency of the considered 
methods of global optimisation not only for the search of global minima of 
multimodal functions, but also for the search of minima of unimodal 
functions of complex geometry, by which all known gradient methods are 
tested (deciding the question of ``how widely the methods can be applied'').

Below we present the results of calculations for the programme implementing 
Halton's and LP-search methods providing an opportunity to specify the 
solution, i.e. to find a point of the global minimum with any predefined 
accuracy. Used as a method of specification is one of the best quasi-Newton 
methods -- the method of Davidon-Fletcher-Powell (DFP) with the most widely 
adopted stopping criterion for smooth problems of the type 
$\|{\grad} f(\bar 
{x}^{k})\| < \varepsilon $.

Authors of the programme (in programming language C++) are EugeneV. Antonij 
and Dmitri B. Poljakov.

Results of calculations for the appropriate test functions are presented in 
Tables 3.1 - 3.7.

\paragraph{\bf "Rozenbrock's Function"} ( -2 $ \le x_1    \le $ 2, -2 $ \le   x_2  
\le $ 2 ) :

\begin{table}[h]
\caption{Rozenbrock's Function}
\begin{tabular}
{|p{50pt}|p{85pt}|p{85pt}|p{85pt}|p{85pt}|}
\hline
Number \par of points N& 
LP-search \par Method& 
LP-search Method \par with specification \par ($\varepsilon $=10$^{ - 6})$& 
Halton's Method& 
Halton's Method \par with specification \par ($\varepsilon $=10$^{ - 6})$ \\
\hline
\par 2000& 
f($\bar {x}^{\ast })$=0.0062603 \par $x_1 ^{\ast } \;\;= 1.0078125$ \par $x_2 ^{\ast }  \;\;= 1.0078125$& 
f($\bar {x}^{\ast })$=0.0000000 \par $x_1 ^{\ast } \;\;= 1.0000000$ \par $x_2 ^{\ast }  \;\;= 1.0000000$& 
f($\bar {x}^{\ast })$=0.0034812 \par $x_1 ^{\ast } \;\;= 1.0214844$ \par $x_2 ^{\ast }  \;\;= 1.0489255$& 
f($\bar {x}^{\ast })$=0.0000000 \par $x_1 ^{\ast } \;\;= 1.0000000$ \par 
$x_2 ^{\ast }  \;\;= 1.0000000$ \\
\hline
\par 8192& 
f($\bar {x}^{\ast })$=0.0043641 \par $x_1 ^{\ast }  \;\;= 0.9638672$ \par $x_2 ^{\ast }  \;\;= 0.9345703$& 
f($\bar {x}^{\ast })$=0.0000000 \par $x_1 ^{\ast }  \;\;= 1.0000000$ \par $x_2 ^{\ast }  \;\;=  1.0000000$& 
f($\bar {x}^{\ast })$=0.0004207 \par $x_1 ^{\ast }  \;\;= 0.9980469$ \par $x_2 ^{\ast }  \;\;=  0.9940558$& 
f($\bar {x}^{\ast })$=0.0000000 \par $x_1 ^{\ast }  \;\;= 0.9999994$ \par 
$x_2 ^{\ast }  \;\;=  0.9999988$ \\
\hline
\par 32767& 
f($\bar {x}^{\ast })$=0.0012168 \par $x_1 ^{\ast }  \;\;= 0.9713135$ \par $x_2 ^{\ast }  \;\;=  0.9454346$& 
f($\bar {x}^{\ast })$=0.0000000 \par $x_1 ^{\ast }  \;\;= 1.0000000$ \par $x_2 ^{\ast }  \;\;=  1.0000001$& 
f($\bar {x}^{\ast })$=0.0004207 \par $x_1 ^{\ast }  \;\;= 0.9980469$ \par $x_2 ^{\ast }  \;\;=  0.9940558$& 
f($\bar {x}^{\ast })$=0.0000000 \par $x_1 ^{\ast }  \;\;= 0.9999994$ \par 
$x_2 ^{\ast }  \;\;=  0.9999988$ \\
\hline
\par 65535& 
f($\bar {x}^{\ast })$=0.0000036 \par $x_1 ^{\ast }  \;\;= 0.9999390$ \par $x_2 ^{\ast }  \;\;= 1.0000610$& 
f($\bar {x}^{\ast })$=0.0000000 \par $x_1 ^{\ast }  \;\;= 1.0000000$ \par $x_2 ^{\ast }  \;\;= 1.0000000$& 
f($\bar {x}^{\ast })$=0.0003650 \par $x_1 ^{\ast }  \;\;= 1.0169067$ \par $x_2 ^{\ast }  \;\;= 1.0332097$& 
f($\bar {x}^{\ast })$=0.0000000 \par $x_1 ^{\ast }  \;\;= 1.0000000$ \par $x_2 ^{\ast }  \;\;= 1.0000000$ \\
\hline
\end{tabular}
\label{tab2}
\end{table}

Solving the given problem of minimization with the help of the MATLAB 6 
system has produced the following results:

\noindent
given $\bar {x}^{0}$ = (-5;-5): $\bar {x}^{\ast }$ = (1.0000000; 
1.0000000); f($\bar {x}^{\ast })$ = 5.6304 $ \cdot $ 10$^{ - 10}$;

\noindent
given $\bar {x}^{0}$ = (5; 5): $\bar {x}^{\ast }$ = (1.0000000; 
1.0000000); f($\bar {x}^{\ast })$ = 5.6197 $ \cdot $ 10$^{ - 10}$;

\noindent
given $\bar {x}^{0}$ = (-5; 5): $\bar {x}^{\ast }$ = (1.0000000; 
1.0000000); f($\bar {x}^{\ast })$ = 8.8830 $ \cdot $ 10$^{ - 10}$.

The result is satisfactory.

\paragraph{\bf "Fletcher-Powell Function"} (-1$ \le x_1  \le $ 1, 0$ \le x_2  \le $ 
2, 0$ \le  x_3  \le $ 2):

\begin{table}[h]
\caption{Fletcher-Powell Function}
\begin{tabular}
{|p{50pt}|p{85pt}|p{85pt}|p{85pt}|p{85pt}|}
\hline
Number of points N& 
LP-search \par Method& 
LP-search Method with specification ($\varepsilon $=10$^{ - 6})$& 
Halton's \par Method& 
Halton's Method with specification \par ($\varepsilon $=10$^{ - 6})$ \\
\hline
2000& 
f($\bar {x}^{\ast })$=0.2782941 \par $x_1 ^{\ast \;\;} =0.9521484$ \par $x_2 ^{\ast \;\;} =0.0410156$ \par $x_3 ^{\ast \;\;} =0.0449219$& 
f($\bar {x}^{\ast })$=0.0000000 \par $x_1 ^{\ast \;\;} =1.0000000$ \par $x_2 ^{\ast \;\;} =0.0000000$ \par $x_3 ^{\ast \;\;} =0.0000000$& 
f($\bar {x}^{\ast })$=0.5961633 \par $x_1 ^{\ast \;\;} =0.9389648$ \par $x_2 ^{\ast \;\;} =0.3063557$ \par $x_3 ^{\ast \;\;} =0.4396800$& 
f($\bar {x}^{\ast })$=0.0000000 \par $x_1 ^{\ast \;\;} =1.0000000$ \par 
$x_2 ^{\ast \;\;} =0.0000000$ \par $x_3 ^{\ast \;\;} =0.0000000$ \\
\hline
8192& 
f($\bar {x}^{\ast })$=0.1721926 \par $x_1 ^{\ast \;\;} =0.9700928$ \par $x_2 ^{\ast \;\;} =0.2556152$ \par $x_3 ^{\ast \;\;} =0.3991699$& 
f($\bar {x}^{\ast })$=0.0000000 \par $x_1 ^{\ast \;\;} =1.0000000$ \par $x_2 ^{\ast \;\;} =0.0000000$ \par $x_3 ^{\ast \;\;} =0.0000000$& 
f($\bar {x}^{\ast })$=0.3818622 \par $x_1 ^{\ast \;\;} =0.9422607$ \par $x_2 ^{\ast \;\;} =0.3713865$ \par $x_3 ^{\ast \;\;} =0.6024960$& 
f($\bar {x}^{\ast })$=0.0000000 \par $x_1 ^{\ast \;\;} =1.0000000$ \par $x_2 ^{\ast \;\;} =0.0000000$ \par $x_3 ^{\ast \;\;} =0.0000000$ \\
\hline
32767& 
f($\bar {x}^{\ast })$=0.0970329 \par $x_1 ^{\ast \;\;} =0.9846497$ \par $x_2 ^{\ast \;\;} =0.1121216$ \par $x_3 ^{\ast \;\;} =0.1549683$& 
f($\bar {x}^{\ast })$=0.0000000 \par $x_1 ^{\ast \;\;} =1.0000000$ \par $x_2 ^{\ast \;\;} =0.0000001$ \par $x_3 ^{\ast \;\;} =0.0000002$& 
f($\bar {x}^{\ast })$=0.0973852 \par $x_1 ^{\ast \;\;} =0.9801636$ \par $x_2 ^{\ast \;\;} =0.0120917$ \par $x_3 ^{\ast \;\;} =0.0433920$& 
f($\bar {x}^{\ast })$=0.0000000 \par $x_1 ^{\ast \;\;} =1.0000000$ \par $x_2 ^{\ast \;\;} =0.0000000$ \par $x_3 ^{\ast \;\;} =0.0000001$ \\
\hline
65535& 
f($\bar {x}^{\ast })$=0.0867985 \par $x_1 ^{\ast \;\;} =0.9713898$ \par $x_2 ^{\ast \;\;} =0.0281677$ \par $x_3 ^{\ast \;\;} =0.0528259$& 
f($\bar {x}^{\ast })$=0.0000000 \par $x_1 ^{\ast \;\;} =1.0000000$ \par $x_2 ^{\ast \;\;} =0.0000000$ \par $x_3 ^{\ast \;\;} =0.0000000$& 
f($\bar {x}^{\ast })$=0.0883128 \par $x_1 ^{\ast \;\;} =0.9718781$ \par $x_2 ^{\ast \;\;} =0.0985283$ \par $x_3 ^{\ast \;\;} =0.1496832$& 
f($\bar {x}^{\ast })$=0.0000000 \par $x_1 ^{\ast \;\;} =1.0000000$ \par $x_2 ^{\ast \;\;} =0.0000000$ \par $x_3 ^{\ast \;\;} =0.0000000$ \\
\hline
\end{tabular}
\label{tab3}
\end{table}

Solving the given problem of minimization with the help of the MATLAB 6 
system has produced the following results:

\noindent
given $\bar {x}^{0}$ = (-10; -10; -10): $\bar {x}^{\ast }$ = (1.0000; 
0.0000; 0.0000); f($\bar {x}^{\ast })$=7.0274 $ \cdot $ 10$^{ - 9}$;

\noindent
given $\bar {x}^{0}$ = (-0.001; 10; 10): $\bar {x}^{\ast }$ = (-0.0017; 
-0.9998; 7.4231); f($\bar {x}^{\ast })$=55.6538.

One can conclude, that in the MATLAB 6 system the problem cannot be solved, 
as there is dependence of the results on the initial conditions of the 
search (which is typical for the local methods of optimisation, or in the 
case of discontinuity of partial derivatives of objective function for some 
variables).

\paragraph{\bf "Powell's Function"} (-1 $ \le  x_i  \le $ 2, i = 1,2,3,4):

\begin{table}[htbp]
\caption{Powell's Function}
\begin{tabular}
{|p{50pt}|p{85pt}|p{85pt}|p{85pt}|p{85pt}|}
\hline
Number of points N& 
LP-search \par Method& 
LP-search Method with specification ($\varepsilon $=10$^{ - 6})$& 
Halton's \par Method& 
Halton's Method with specification \par ($\varepsilon $=10$^{ - 6})$ \\
\hline
2000 \par & 
f($\bar {x}^{\ast })$=0.6550502 \par $x_1 ^{\ast \;\;} =- 0.7890625$ \par $x_2 ^{\ast \;\;} = 0.0546875$ \par $x_3 ^{\ast \;\;} =- 0.3671875$ \par $x_4 ^{\ast \;\;} = - 0.4140625$& 
f($\bar {x}^{\ast })$=0.0000000 \par $x_1 ^{\ast \;\;} =- 0.0001419$ \par $x_2 ^{\ast \;\;} =0.0000142$ \par $x_3 ^{\ast \;\;} =- 0.0001390$ \par $x_4 ^{\ast \;\;} = - 0.0001389$& 
f($\bar {x}^{\ast })$=0.7399679 \par $x_1 ^{\ast \;\;} =0.3549805$ \par
$x_2 \quad =0.0397805$ \par $x_3 ^{\ast \;\;} = - 0.0572800$ \par $x_4 ^{\ast \;\;} = 0.0020825$& 
f($\bar {x}^{\ast })$=0.0000000 \par $x_1 ^{\ast \;\;} =0.0003630$ \par $x_2 ^{\ast \;\;} = - 0.0000363$ \par $x_3 ^{\ast \;\;} =0.0003281$ \par $x_4 ^{\ast \;\;} =0.0003280$ \\
\hline
8192. \par 32767& 
f($\bar {x}^{\ast })$=0.1618187 \par $x_1 ^{\ast \;\;} =- 0.0518799$ \par $x_2 ^{\ast \;\;} = 0.0360107$ \par $x_3 ^{\ast \;\;} =- 0.1236572$ \par $x_4 ^{\ast \;\;} = - 0.0137939$& 
f($\bar {x}^{\ast })$=0.0000000 \par $x_1 ^{\ast \;\;} =0.0000816$ \par $x_2 ^{\ast \;\;} = - 0.0000082$ \par $x_3 ^{\ast \;\;} =0.0002766$ \par $x_4 ^{\ast \;\;} = 0.0002766$& 
f($\bar {x}^{\ast })$=0.0711034 \par $x_1 ^{\ast \;\;} = - 0.1427002$ \par $x_2 ^{\ast } =0.0269776$ \par $x_3 ^{\ast \;\;} = - 0.0273280$ \par $x_4 ^{\ast \;\;} = 0.0599155$& 
f($\bar {x}^{\ast })$=0.0000000 \par $x_1 ^{\ast \;\;} = - 0.0001949$ \par $x_2 ^{\ast \;\;} =0.0000195$ \par $x_3 ^{\ast \;\;} = - 0.0000995$ \par $x_4 ^{\ast \;\;} = - 0.0000995$ \\
\hline
65535& 
f($\bar {x}^{\ast })$=0.0625249 \par $x_1 ^{\ast \;\;} =0.2141266$ \par $x_2 ^{\ast \;\;} = - 0.0077057$ \par $x_3 ^{\ast \;\;} =0.0576630$ \par $x_4 ^{\ast \;\;} =0.1507721$& 
f($\bar {x}^{\ast })$=0.0000000 \par $x_1 ^{\ast \;\;} =0.0000578$ \par $x_2 ^{\ast \;\;} = - 0.0000058$ \par $x_3 ^{\ast \;\;} =- 0.0002014$ \par $x_4 ^{\ast \;\;} = - 0.0002014$& 
f($\bar {x}^{\ast })$=0.0336186 \par $x_1 ^{\ast \;\;} = - 0.0893707$ \par $x_2 ^{\ast } = - 0.0030991$ \par $x_3 ^{\ast \;\;} =0.0338045$ \par $x_4 ^{\ast \;\;} = 0.0803152$& 
f($\bar {x}^{\ast })$=0.0000000 \par $x_1 ^{\ast \;\;} =0.0028412$ \par $x_2 ^{\ast \;\;} = - 0.0002841$ \par $x_3 ^{\ast \;\;} =0.0014357$ \par $x_4 ^{\ast \;\;} =0.0014357$ \\
\hline
\end{tabular}
\label{tab4}
\end{table}

Solving the given problem of minimization with the help of the MATLAB 6 
system has produced the following results:

1) given $\bar {x}^{0}$ = (-5;-5;-5;-5): $\bar {x}^{\ast }$ = 
(0.0002799; 0.0000280; -0.0000742; -0.0000743); 

\noindent
f($\bar {x}^{\ast })$ = 5.0413 $ \cdot $ 10$^{ - 14}$; 

2) given $\bar {x}^{0}$ = (-5; 5;-5; 5): $\bar {x}^{\ast }$ = 
(-0.0001917; 0.0000192; 0.0006634; 0.0006633);

\noindent
f($\bar {x}^{\ast })$ = 8.3379 $ \cdot $ 10$^{ - 12}$.

The result is satisfactory.

\paragraph{\bf "Wood's Function"} (0 $ \le  x_i   \le $ 3, i = 1,2,3,4):

\begin{table}[htbp]
\caption{Wood's Function}
\begin{tabular}
{|p{50pt}|p{85pt}|p{85pt}|p{85pt}|p{85pt}|}
\hline
Number of points N& 
{LP-search \par Method} & 
{LP-search Method with specification ($\varepsilon $=10$^{ - 6})$} & 
Halton's \par Method& 
Halton's Method with specification \par ($\varepsilon $=10$^{ - 6})$ \\
\hline
2000 \par & 
f($\bar {x}^{\ast })$=2.6487434 \par $x_1 ^{\ast \;\;} =0.8862305$ \par $x_2 ^{\ast \;\;} = 0.7163086$ \par $x_3 ^{\ast \;\;} =1.3051758$ \par $x_4 ^{\ast \;\;} = 1.5981445$& 
{f($\bar {x}^{\ast })$=0.0000000 \par $x_1 ^{\ast \;\;} =1.0000000$ \par $x_2 ^{\ast \;\;} =1.0000000$ \par $x_3 ^{\ast \;\;} =1.0000000$ \par $x_4 ^{\ast \;\;} =1.0000000$} & 
{f($\bar {x}^{\ast })$=3.3474517 \par $x_1 ^{\ast \;\;} =1.1176758$ \par $x_2 ^{\ast \;\;} =1.3525377$ \par $x_3 ^{\ast \;\;} =0.8851200$ \par $x_4 ^{\ast \;\;} =0.6272387$} & 
f($\bar {x}^{\ast })$=0.0000000 \par $x_1 ^{\ast \;\;} =1.0000000$ \par $x_2 ^{\ast \;\;} =1.0000000$ \par $x_3 ^{\ast \;\;} =1.0000000$ \par $x_4 ^{\ast \;\;} =1.0000000$ \\
\hline
8192& 
f($\bar {x}^{\ast })$=1.7646029 \par $x_1 ^{\ast \;\;} =0.09338379$ \par $x_2 ^{\ast \;\;} = 0.0186768$ \par $x_3 ^{\ast \;\;} =1.4117432$ \par $x_4 ^{\ast \;\;} =1.9288330$& 
{f($\bar {x}^{\ast })$=0.0000000 \par $x_1 ^{\ast \;\;} =1.0000000$ \par $x_2 ^{\ast \;\;} =1.0000000$ \par $x_3 ^{\ast \;\;} =1.0000000$ \par $x_4 ^{\ast \;\;} =1.0000000$} & 
{f($\bar {x}^{\ast })$=0.8959908 \par $x_1 ^{\ast \;\;} =0.8909912$ \par $x_2 ^{\ast } =0.7892090$ \par $x_3 ^{\ast \;\;} =1.0955520$ \par $x_4 ^{\ast \;\;} =1.1084667$} & 
f($\bar {x}^{\ast })$=0.0000000 \par $x_1 ^{\ast \;\;} =1.0000000$ \par $x_2 ^{\ast \;\;} =1.0000000$ \par $x_3 ^{\ast \;\;} =1.0000000$ \par $x_4 ^{\ast \;\;} =1.0000000$ \\
\hline
32767& 
f($\bar {x}^{\ast })$=1.4175885 \par $x_1 ^{\ast \;\;} =1.1044006$ \par $x_2 ^{\ast \;\;} =1.1336975$ \par $x_3 ^{\ast \;\;} =0.9775085$ \par $x_4 ^{\ast \;\;} =1.0236511$& 
{f($\bar {x}^{\ast })$=0.0000000 \par $x_1 ^{\ast \;\;} =1.0000000$ \par $x_2 ^{\ast \;\;} =1.0000000$ \par $x_3 ^{\ast \;\;} =1.0000000$ \par $x_4 ^{\ast \;\;} =1.0000000$} & 
{f($\bar {x}^{\ast })$=0.8959908 \par $x_1 ^{\ast \;\;} =0.8909912$ \par $x_2 ^{\ast } =0.7892090$ \par $x_3 ^{\ast \;\;} =1.0955520$ \par $x_4 ^{\ast \;\;} =1.1084667$} & 
f($\bar {x}^{\ast })$=0.0000000 \par $x_1 ^{\ast \;\;} =1.0000000$ \par $x_2 ^{\ast \;\;} =1.0000000$ \par $x_3 ^{\ast \;\;} =1.0000000$ \par $x_4 ^{\ast \;\;} =1.0000000$ \\
\hline
65535& 
f($\bar {x}^{\ast })$=0.1698878 \par $x_1 ^{\ast \;\;} =0.7766876$ \par $x_2 ^{\ast \;\;} =0.6017303$ \par $x_3 ^{\ast \;\;} =1.1534271$ \par $x_4 ^{\ast \;\;} =1.3391876$& 
{f($\bar {x}^{\ast })$=0.0000000 \par $x_1 ^{\ast \;\;} =1.0000000$ \par $x_2 ^{\ast \;\;} =1.0000000$ \par $x_3 ^{\ast \;\;} =1.0000000$ \par $x_4 ^{\ast \;\;} =1.0000000$} & 
{f($\bar {x}^{\ast })$=0.8959908 \par $x_1 ^{\ast \;\;} =0.8909912$ \par $x_2 ^{\ast } =0.7892090$ \par $x_3 ^{\ast \;\;} =1.0955520$ \par $x_4 ^{\ast \;\;} =1.1084667$} & 
f($\bar {x}^{\ast })$=0.0000000 \par $x_1 ^{\ast \;\;} =1.0000000$ \par $x_2 ^{\ast \;\;} =1.0000000$ \par $x_3 ^{\ast \;\;} =1.0000000$ \par $x_4 ^{\ast \;\;} =1.0000000$ \\
\hline
\end{tabular}
\label{tab5}
\end{table}

Solving the given problem of minimization with the help of the MATLAB 6 
system has produced the following results:

\noindent
given $\bar {x}^{0}$ = (-5; 4; -3; 2): $\bar {x}^{\ast }$ = (-0.2067; 
0.0397; 1.4070; 1.9820); f($\bar {x}^{\ast })$=2.0051;

\noindent
given $\bar {x}^{0}$ = (-2;-2;-3;-3): $\bar {x}^{\ast }$ = (0.6295; 
0.4334;-1.3062; 1.6882); f($\bar {x}^{\ast })$=5.9277;

\noindent
given $\bar {x}^{0}$ = ( 0; 0; 0; 0): $\bar {x}^{\ast }$ = (1.0000; 
1.0000; 1.0000; 1.0000); f($\bar {x}^{\ast })$=6.2146 $ \cdot $ 10$^{ - 
8}$.

One can conclude, that in system MATLAB 6 the problem cannot be solved, as 
there is dependence of the results on the initial conditions of the search 
(which is typical for the local methods of optimisation).

\paragraph{\bf Function No. 5} (-3 $ \le x_1 \le $ 1, -1 $ \le x_2  \le $ 3):

\begin{table}[htbp]
\caption{Function No. 5}
\begin{tabular}
{|p{50pt}|p{85pt}|p{85pt}|p{85pt}|p{85pt}|}
\hline
Number of points N& 
{LP-search \par Method} & 
{LP-search Method with specification ($\varepsilon $=10$^{ - 6})$} & 
Halton's \par Method& 
Halton's Method with specification \par ($\varepsilon $=10$^{ - 6})$ \\
\hline
2000& 
f($\bar {x}^{\ast })$=-2.000000 $x_1 ^{\ast }$= 0.000000 \par $x_2 ^{\ast \;\;} =0.000000$& 
{f($\bar {x}^{\ast })$=-2.000000 \par $x_1 ^{\ast \;\;} =0.000000$ \par $x_2 ^{\ast \;\;} =0.000000$} & 
{f($\bar {x}^{\ast })$=-1.817494 \par $x_1 ^{\ast \;\;} =- 0.363281$ \par $x_2 ^{\ast \;\;} =- 0.010517$} & 
f($\bar {x}^{\ast })$=-1.878901 \par $x_1 ^{\ast \;\;} =- 0.346924$ \par $x_2 ^{\ast \;\;} =- 0.000000$ \\
\hline
8192& 
f($\bar {x}^{\ast })$=-2.000000 \par $x_1 ^{\ast \;\;} =0.000000$ \par $x_2 ^{\ast \;\;} = 0.000000$& 
{f($\bar {x}^{\ast })$=-2.000000 \par $x_1 ^{\ast \;\;} =0.000000$ \par $x_2 ^{\ast \;\;} =0.000000$} & 
{f($\bar {x}^{\ast })$=-1.860544 \par $x_1 ^{\ast \;\;} =0.008789$ \par $x_2 ^{\ast \;\;} = - 0.340954$} & 
f($\bar {x}^{\ast })$=-1.878901 \par $x_1 ^{\ast \;\;} =- 0.000000$ \par $x_2 ^{\ast \;\;} =- 0.346924$ \\
\hline
32767& 
f($\bar {x}^{\ast })$=-2.000000 \par $x_1 ^{\ast \;\;} =0.000000$ \par $x_2 ^{\ast \;\;} =0.000000$& 
{f($\bar {x}^{\ast })$=-2.000000 \par $x_1 ^{\ast \;\;} =0.000000$ \par $x_2 ^{\ast \;\;} =0.000000$} & 
{f($\bar {x}^{\ast })$=-1.963421 \par $x_1 ^{\ast \;\;} =- 0.012817$ \par $x_2 ^{\ast \;\;} =- 0.007807$} & 
f($\bar {x}^{\ast })$=-2.000000 \par $x_1 ^{\ast \;\;} =0.000000$ \par $x_2 ^{\ast \;\;} =- 0.000000$ \\
\hline
65535& 
f($\bar {x}^{\ast })$=-2.000000 \par $x_1 ^{\ast \;\;} =0.000000$ \par $x_2 ^{\ast \;\;} =0.000000$& 
{f($\bar {x}^{\ast })$=-2.000000 \par $x_1 ^{\ast \;\;} =0.000000$ \par $x_2 ^{\ast \;\;} =0.000000$} & 
{f($\bar {x}^{\ast })$=-1.996424 \par $x_1 ^{\ast }$ = 0.004578 \par $x_2 
^{\ast \;\;} =0.000999$} & 
f($\bar {x}^{\ast })$=-2.000000 \par $x_1 ^{\ast \;\;} =- 0.000000$ \par $x_2 ^{\ast \;\;} =- 0.000000$ \\
\hline
\end{tabular}
\label{tab6}
\end{table}

Solving the given problem of minimization with the help of the MATLAB 6 
system has produced the following results:

\noindent
given $\bar {x}^{0}$ = (-3;-3): $\bar {x}^{\ast }$ = (-3.1219; -3.1219); 
f($\bar {x}^{\ast })$= 17.6169;

\noindent
given $\bar {x}^{0}$ = (-1;-1): $\bar {x}^{\ast }$ = (-1.0408; -1.0408); 
f($\bar {x}^{\ast })$= 0.1798;

\noindent
given $\bar {x}^{0}$ = ( 0.5; -0.5): $\bar {x}^{\ast }$ = (-0.3469; 
0.3469); f($\bar {x}^{\ast })$= -1.7578;

\noindent
given $\bar {x}^{0}$ = ( -0.5; 0.5): $\bar {x}^{\ast }$ = (0.3469; 
-0.3469); f($\bar {x}^{\ast })$= -1.7578;

\noindent
given $\bar {x}^{0}$ = (0.2; 0.2): $\bar {x}^{\ast }$ = (0.3469; 
0.3469); f($\bar {x}^{\ast })$= -1.7578;

\noindent
given $\bar {x}^{0}$ = (0.1; 0.1): $\bar {x}^{\ast }$ = (0.0000113; 
0.0000113); f($\bar {x}^{\ast })$= -2.000000;

\noindent
given $\bar {x}^{0}$ = (0.05; 0.05): $\bar {x}^{\ast }$ = (-0.000026; 
-0.000026); f($\bar {x}^{\ast })$= -2.000000;

One can conclude, that in the MATLAB 6 system the problem cannot be solved, 
as there is dependence of the results on the initial conditions of the 
search (the system finds local minima of the function within the vicinity of 
a given starting point).

\paragraph{\bf ``Himmelblau's function''} (2.002 $ \le x_i  \le $ 9.998; i = 
1,2,\ldots ,10):

\begin{table}[htbp]
\caption{Himmelblau's function}
\begin{tabular}
{|p{50pt}|p{85pt}|p{85pt}|p{85pt}|p{85pt}|}
\hline
Amount of points N& 
LP-search \par Method& 
LP-search Method with specification ($\varepsilon $=10$^{ - 6})$& 
Halton's \par Method& 
Halton's Method with specification \par ($\varepsilon $=10$^{ - 6})$ \\
\hline
2000 \par & 
f($\bar {x}^{\ast })$=-21.28223 \par $x_1 ^{\ast \;\;} =6.132746$ \par $x_2 ^{\ast \;\;} =6.398238$ \par $x_3 ^{\ast \;\;} =8.006809$ \par $x_4 ^{\ast \;\;} =8.537793$ \par $x_5 ^{\ast \;\;} =9.552910$ \par $x_6 ^{\ast \;\;} =7.897488$ \par $x_7 ^{\ast \;\;} =8.537793$ \par $x_8 ^{\ast \;\;} =5.757934$ \par $x_9 ^{\ast \;\;} =9.131246$ \par $x_{10} ^{\ast } =8.428473$ \par & 
f($\bar {x}^{\ast })$=-45.77846 \par $x_1 ^{\ast \;\;} =9.350265$ \par $x_2 ^{\ast \;\;} =9.350265$ \par $x_3 ^{\ast \;\;} =9.350265$ \par $x_4 ^{\ast \;\;} =9.350265$ \par $x_5 ^{\ast \;\;} =9.350265$ \par $x_6 ^{\ast \;\;} =9.350265$ \par $x_7 ^{\ast \;\;} =9.350265$ \par $x_8 ^{\ast \;\;} =9.350265$ \par $x_9 ^{\ast \;\;} =9.350265$ \par $x_{10} ^{\ast } =9.350265$& 
f($\bar {x}^{\ast })$=-24.99797 \par $x_1 ^{\ast \;\;} =9.798881$ \par $x_2 ^{\ast \;\;} =9.398390$ \par $x_3 ^{\ast \;\;} =9.184327$ \par $x_4 ^{\ast \;\;} =6.934143$ \par $x_5 ^{\ast \;\;} =9.139471$ \par $x_6 ^{\ast \;\;} =7.850690$ \par $x_7 ^{\ast \;\;} =7.620776$ \par $x_8 ^{\ast \;\;} =8.295979$ \par $x_9 ^{\ast \;\;} =6.764638$ \par $x_{10} ^{\ast } =7.155517$ \par & 
f($\bar {x}^{\ast })$=-45.77846 \par $x_1 ^{\ast \;\;} =9.350265$ \par $x_2 ^{\ast \;\;} =9.350265$ \par $x_3 ^{\ast \;\;} =9.350265$ \par $x_4 ^{\ast \;\;} =9.350265$ \par $x_5 ^{\ast \;\;} =9.350265$ \par $x_6 ^{\ast \;\;} =9.350265$ \par $x_7 ^{\ast \;\;} =9.350265$ \par $x_8 ^{\ast \;\;} =9.350265$ \par $x_9 ^{\ast \;\;} =9.350265$ \par $x_{10} ^{\ast } =9.350265$ \\
\hline
65535& 
f($\bar {x}^{\ast })$=-37.85620 \par $x_1 ^{\ast \;\;} =8.814022$ \par $x_2 ^{\ast \;\;} =9.325485$ \par $x_3 ^{\ast \;\;} =9.040471$ \par $x_4 ^{\ast \;\;} =8.489965$ \par $x_5 ^{\ast \;\;} =8.763266$ \par $x_6 ^{\ast \;\;} =7.385049$ \par $x_7 ^{\ast \;\;} =8.409927$ \par $x_8 ^{\ast \;\;} =8.878443$ \par $x_9 ^{\ast \;\;} =9.042423$ \par $x_{10} ^{\ast } =8.958481$& 
f($\bar {x}^{\ast })$=-45.77846 \par $x_1 ^{\ast \;\;} =9.350265$ \par $x_2 ^{\ast \;\;} =9.350265$ \par $x_3 ^{\ast \;\;} =9.350265$ \par $x_4 ^{\ast \;\;} =9.350265$ \par $x_5 ^{\ast \;\;} =9.350265$ \par $x_6 ^{\ast \;\;} =9.350265$ \par $x_7 ^{\ast \;\;} =9.350265$ \par $x_8 ^{\ast \;\;} =9.350265$ \par $x_9 ^{\ast \;\;} =9.350265$ \par $x_{10} ^{\ast } =9.350265$& 
f($\bar {x}^{\ast })$=-35.07950 \par $x_1 ^{\ast \;\;} =9.319751$ \par $x_2 ^{\ast \;\;} =8.630329$ \par $x_3 ^{\ast \;\;} =7.693310$ \par $x_4 ^{\ast \;\;} =9.382510$ \par $x_5 ^{\ast \;\;} =8.223152$ \par $x_6 ^{\ast \;\;} =8.532420$ \par $x_7 ^{\ast \;\;} =7.939859$ \par $x_8 ^{\ast \;\;} =8.213517$ \par $x_9 ^{\ast \;\;} =8.880209$ \par $x_{10} ^{\ast } =8.766924$& 
f($\bar {x}^{\ast })$=-45.77846 \par $x_1 ^{\ast \;\;} =9.350265$ \par $x_2 ^{\ast \;\;} =9.350265$ \par $x_3 ^{\ast \;\;} =9.350265$ \par $x_4 ^{\ast \;\;} =9.350265$ \par $x_5 ^{\ast \;\;} =9.350265$ \par $x_6 ^{\ast \;\;} =9.350265$ \par $x_7 ^{\ast \;\;} =9.350265$ \par $x_8 ^{\ast \;\;} =9.350265$ \par $x_9 ^{\ast \;\;} =9.350265$ \par $x_{10} ^{\ast } =9.350265$ \\
\hline
\end{tabular}
\label{tab7}
\end{table}

Solving the given problem of minimization with the help of the MATLAB 6 
system has produced the following results:

\noindent
given $\bar {x}^{0}$ = (3;3;\ldots ;3): the task does not have solutions;

\noindent
given $\bar {x}^{0}$ = (5;5;\ldots ;5): 

$\bar {x}^{\ast \;\;} = $(7.6458; 9.0869; 9.3527; 9.3878; 9.3332; 9.2646; 
9.2418; 9.2781; 9.3716; 9.3649); f($\bar {x}^{\ast })$= -42.5434;

\noindent
given $\bar {x}^{0}$ = (8; 3; 4; 5; 6; 7; 3; 5; 6; 7): 

$\bar {x}^{\ast \;\;} = $(9.1030; 8.6182; 5.4010; 9.3945; 8.6888; 9.3815; 
9.5639; 9.2224; 9.2201; 9.4322); f($\bar {x}^{\ast })$= -35.1252;

\noindent
given $\bar {x}^{0}$ = (3; 4; 5; 6; 7; 8; 3; 4; 5; 6): the task does not 
have solutions;

\noindent
given $\bar {x}^{0}$ = (9; 8; 7; 6; 5; 9; 8; 7; 6; 5): 

$\bar {x}^{\ast \;\;} = $(9.2172; 9.1307; 9.3271; 9.2110; 9.2809; 4.5076; 
9.2679; 9.2796; 8.6958; 8.4667);

\noindent
f($\bar {x}^{\ast })$= -32.6965.

Here we conclude, that in the MATLAB 6 system the problem cannot be solved.

For the test function  7 the solution was found by the method of LP-search 
and by Halton's method without the procedure of specification, for N = 
65535. The following results were obtained:

\noindent
method of LP-search: f($\bar {x}^{\ast })$ = 27845.37; $\bar {x}^{\ast 
}$ = (1.49955; 6.12384); 

Halton's method: f($\bar {x}^{\ast })$ = 27845.02; $\bar {x}^{\ast }$ = 
(1.50398; 6.14608).

\noindent
which are in good agreement with the known solution of the given problem.

The results of the calculations show extremely high efficiency of Halton's 
and LP-search methods for solving a problem of numerical global optimisation 
of a function of many variables.

Also, the advantage of the methods is that neither the differentiability of 
the objective function, nor even its analytical definability is required for 
their application. It is only enough to have an opportunity to calculate the 
values of the function at any arbitrary points in its domain, which provides 
for a wide opportunity of application of the methods and the ability to 
solve the majority of practical problems.

At the same time the advantages of gradients methods should not be 
underestimated, the main one being the high speed of convergence to a point 
of a minimum. It appears to be expedient to use the better qualities of 
direct and gradient methods and to use them in combination, which can 
substantially reduce the volume of the calculations necessary for the search 
for the global minimum with a specified accuracy. This conclusion is also 
supported by the results of the numerical experiments presented above.

\section{Final remarks}

In summary we shall consider one more problem of numerical optimisation -- 
the problems of the large dimensions. It is known [4], that in the case when 
the dimension \textbf{n} of the argument of the objective function becomes 
very large, application of gradient procedures becomes impossible or causes 
significant technical difficulties.

If the dimension of the objective function \textbf{n} $>$ 51, it is then 
possible to take advantage of generalized LP$_{\tau }$ -sequence calculated 
with the same formulas as for the case $1 \le  j < \infty $. 
Unfortunately, there are no simple-enough formulas for the calculation of 
the directing points, as was already mentioned above. Therefore it is 
possible to take advantage of the already developed table of numbers $V_j 
^s$ for \textbf{n} $ \le $ 51, and use the usual pseudorandom numbers 
\textbf{$\gamma $ }for the missing coordinates of quasi-random 
points\textbf{, }for example, Q$_{i}$ = (q$_{1}^{i}$, q$_{2}^{i}$, 
\ldots , q$_{51}^{i}$, $\gamma _{52}^{i}$, $\gamma _{53}^{i}$, 
\ldots , $\gamma _{n}^{i})$. As various coordinates of points Q$_{i}$ 
are disparate [11] (coordinates with smaller indexes are better 
distributed), it is useful to assign indexes to the variables of the 
objective function in a way, when the most essential coordinates have 
smaller indexes, for which q$_{j}^{i}$ is calculated, and for all other 
coordinates - $\gamma _{j}^{i}$. Such way of calculation can speed up 
the convergence, in comparison with the calculation only by random points 
$\gamma $.

The second approach to the problem of global optimisation in case of a large 
dimension of the objective function consists of the use of the generalized 
Halton's sequence {\{}P$_{i}${\}}, i = 1, 2, \ldots , the points of which 
have coordinates P$_{i}$ = (p$_{r_1} $(i), p$_{r_2} $(i), \ldots , 
p$_{r_n} $(i),\ldots ), where r$_{1} < r_{2} < r_{3} < \ldots < 
r_{n} < \ldots$ -- is a sequence of all prime numbers. For the calculation 
of the points P$_{i}$, in practice, it is possible to set a big enough table 
of prime numbers or to program an algorithm of their finding, for example, 
the Eratosthenes' Sieve method. The method consists of successively deleting 
all the composite numbers from a sequence of natural numbers. For example, 
let it be required to find all the prime numbers between 1 and 30. For this 
purpose one needs to write down all natural numbers from 1 to 30 in 
ascending order. First number is 1 -- not a prime number, therefore it is 
deleted. The following is number 2 -- prime, it is kept, and every second 
number after 2, i.e. 4, 6, 8, \ldots is deleted. The following prime number 
3 is kept, and every third number, after 3, is deleted. The following prime 
number 5 is kept, and every fifth number, after 5, is deleted (numbers 
already deleted are also taken into account) etc. As a result all prime 
numbers, smaller than 30, are obtained (the kept undeleted numbers): 2, 3, 
5, 7, 11, 13, 17, 19, 23, 29. If it is necessary to find all prime numbers 
not exceeding $N$ than one has to use the procedure described above up to the 
greatest prime number $p$, not exceeding $\sqrt N $. For example, if it is 
necessary to create a table of prime numbers not exceeding 1000, it is 
necessary to finish the deletion procedure on the number 31 inclusive. There 
are now printed tables of prime numbers for up to 12 million, i.e. it is 
possible to obtain sequences of Halton's points practically for any 
dimension.

\section*{References}

\begin{enumerate}

\item
G.V. Reklaitis, A.Ravindran, K.M.Ragsdell. Engineering Optimisation (Methods 
and Applications). Ò. 1,2 - M.: the World, 1986.
\item
Dadid M.Himmelblau. Applied Nonlinear Programming. - M.: the World, 1975.
\item
Brian D.Bunday Basic Optimisation Methods. - M.: Radio and communication 
(connection), 1988.
\item
Philip E.Gill, Walter Murray, Margaret H.Wright. Practical optimisation. - 
M.: the World, 1985.
\item
Douglass J.Wilde. Optimum seeking methods. M.: Science. 1967.
\item
Terry E.Shoup. A Practical Guide to Computer Methods for Engineers. M.: the 
World, 1982.
\item
Polak E. Computational Methods in Optimisation. M.: the World, 1974.
\item
L.A.Rastrigin. Statistical methods of the search.-M.: Science, 1968.
\item
L.A.Rastrigin. Systems of the extreme control.-M.:PhysMathPub, 1974.
\item
I.M.Sobol'. Multivariate integrals and a method of Monte Carlo. DAN, 114, 
1957.
\item
I.M.Sobol'. Multivariate quadrature formulas and Haar's functions.- M.: 
Science, 1969.
\item
I.M.Sobol'. Numerical methods of Monte Carlo.-M.: Science, 1973
\item
I.M.Sobol', R.B.Statnikov. LP-search and problems of optimum design. 
Problems of random search. Riga, 1972.
\item
I.M.Sobol', R.B.Statnikov. The Choice of optimal parameters in
problems with many criterions.-M.: Science, 1981.
\item
J.H.Halton, D.C.Handscomb,  A method for increasing the efficiency of
Monte Carlo integrations. J.Assoc.Comput.Machinery, 1957, 4, n.3, 329-340.
\item
J.H.Halton,  A retrospective and prospective survey of the Monte Carlo
method. SIAM Rev., 1970, 12, n.1, 1-63.
\item
V.G.Pocket. Mathematical programming.-M.: Science, 1975.
\item
Vasil'ev F.P. Numerical methods of the solving of the extreme problems.-M.: 
Science, 1980.
\item
A.G.Suharev, A.V.Timohov, V.V.Fedorov. A rate of methods of 
optimisation.-M.: Science, 1986.
\item
N.N.Moiseev, U.P Ivanilov, E.M.Stoljarova. Methods of optimisation 
.-M.:Science, 1978.
\item
M.M.Medynsky. Numerical methods of optimisation in problems of design of 
systems of equipment flying devices. M.: MAI. 1993. 
\end{enumerate}

\end{document}